\documentclass[10pt]{article}

\usepackage[draft]{graphics}
\usepackage{graphicx}          	 
\usepackage{bm}                 
\usepackage{amsmath}             
\usepackage{amssymb}
\usepackage{amsfonts}             
\usepackage{verbatim}           
\usepackage{amsthm}             

\usepackage{mathtools}
\usepackage{relsize}

\usepackage{makeidx}

\theoremstyle{plain}             

\newtheorem{theorem}{Theorem}[section]

\theoremstyle{definition}
\newtheorem{example}[theorem]{Example}

\newtheorem{remark}[theorem]{Remark}

\makeatletter
\@addtoreset{equation}{chapter}

\makeatother

\def\eqref#1{(\ref{#1})}
\def\dsp{\displaystyle}
\def\Frac#1#2{\frac
{
 {\raise.6ex
 \hbox{$\displaystyle#1$}}
}
{
 {\lower.6ex
 \hbox{$\displaystyle#2$}}
 }
}

\numberwithin{equation}{section}

\def\bigOxe{\sqcup \kern-2.3mm \sqcap}

\def\eoexample{{\unskip\nobreak\hfil\penalty50	
\hskip2em\hbox{}\nobreak\hfil$\diamondsuit$
\parfillskip=0pt\finalhyphendemerits=0\medbreak}}

\def\eoremark{{\unskip\nobreak\hfil\penalty50	
\hskip2em\hbox{}\nobreak\hfil$\triangle$
\parfillskip=0pt\finalhyphendemerits=0\medbreak}}

\def\eps{\varepsilon}

\makeatother

\def\dsp{\displaystyle}
\def\Frac#1#2{\frac
{
 {\raise.6ex
 \hbox{$\displaystyle#1$}}
}
{
 {\lower.6ex
 \hbox{$\displaystyle#2$}}
 }
}

\def\eps{{\varepsilon}}

\input epsf

\def\CHFs#1#2#3{
{}_1F_1\left({a};{c};{z}\right)
}

\def\bigO{{\cal O}}
\def\calC{{{\cal C}}}

\def\Z z_+{{\mathbb Z^+}}      

\def\phase{{\rm ph}}

\def\tfrac#1#2{{{\lower.6ex
\hbox{$\scriptstyle#1$}}\over 
{\raise.7ex
\hbox{$\scriptstyle#2$}}}}

\begin{document}
 \title{Asymptotic expansions of Jacobi polynomials and of  the nodes and weights of Gauss-Jacobi quadrature 
for large degree and parameters in terms of elementary functions
}

\author{
A. Gil\\
Departamento de Matem\'atica Aplicada y CC. de la Computaci\'on.\\
ETSI Caminos. Universidad de Cantabria. 39005-Santander, Spain.\\
 \and
J. Segura\\
        Departamento de Matem\'aticas, Estadistica y 
        Computaci\'on,\\
        Universidad de Cantabria, 39005 Santander, Spain.\\
\and
N. M. Temme\\
IAA, 1825 BD 25, Alkmaar, The Netherlands.\footnote{Former address: Centrum Wiskunde \& Informatica (CWI), 
        Science Park 123, 1098 XG Amsterdam,  The Netherlands}\\
}


\maketitle
\begin{abstract}
Asymptotic approximations of Jacobi polynomials  are given in terms of elementary functions for large degree $n$ and parameters $\alpha$ and $\beta$.  From these new results, asymptotic expansions of the zeros are derived and methods are given to obtain the coefficients in the expansions. These approximations can be used as initial values in iterative methods for computing the nodes
of Gauss--Jacobi quadrature for  large degree and parameters. The performance of
 the asymptotic approximations for computing the nodes and weights of these Gaussian quadratures is
illustrated with numerical examples.
\end{abstract}

\section{Introduction}\label{sec:Intro}
This paper is a further exploration in our research on Gauss quadrature for the classical orthogonal polynomials; earlier publications are \cite{Gil:2018:FGH}, \cite{Gil:2018:GHL},  \cite{Gil:2018:AJP},  \cite{Gil:2019:NIG}. Other recent relevant papers on this topic are \cite{Bogaert:2014:IFC}, \cite{Hale:2013:FAC}, \cite{Town:2016:IMA}. 

When we assume that the degree $n$ and the two parameters  $\alpha$ and $\beta$ of the Jacobi polynomial $P_n^{(\alpha,\beta)}(x)$ are large, and we consider the variable $x$ as a parameter that causes nonuniform behavior of the polynomial, it can be expected that, for a detailed and optimal description of the asymptotic approximation, we need a function of three variables. Candidates for this are the Gegenbauer and  the Laguerre polynomial. The Gegenbauer polynomial can be used when the ratio $\alpha/\beta$  does not tend to zero or to infinity. When it does, the Laguerre polynomial is the best option. 

It is possible to transform an integral of $P_n^{(\alpha,\beta)}(x)$ into an integral resembling one of the Gegenbauer or the Laguerre polynomial (and similar when we are working with differential equations). From a theoretical point of view this may be of interest, however,  for practical purposes, when using the results for Gauss quadrature, the transformations and  the coefficients in the expansions become rather complicated. In addition, computing the approximants, that is, large degree polynomials with large additional parameter and a variable in domains where nonuniform behavior of these polynomials may happen, gives an extra nontrivial complication.

Even when we use the Bessel functions or Hermite polynomials as approximants, these complications are still quite relevant. For this reason we consider in this paper expansions in terms of elementary functions, and we will see that to evaluate a certain number of coefficients already gives quite complicated expressions.

For large values of $\beta$ with fixed degree $n$ we have quite simple results derived in \cite{Gil:2018:AJP}, which paper is inspired by \cite{Dimitrov:2016:ABJ}. Large-degree results valid near $x=1$ are given in \cite[\S28.4]{Temme:2015:AMI}, and for the case that $\beta$ is large as well we refer to  \cite[\S28.4.1]{Temme:2015:AMI}.

\section{Several asymptotic phenomena}\label{sec:phen}

To describe the behavior of the Jacobi polynomial for large degree and parameters $\alpha$ and $\beta$, with $x\in[-1,1]$, it is instructive to consider the differential equation of the function
\begin{equation}\label{eq:Intro01}
 W(x)=(1-x)^{\frac12(\alpha+1)}(1+x)^{\frac12(\beta+1)}P_n^{(\alpha,\beta)}( x).
\end{equation}
By using the Liouville-Green transformations as described in \cite{Olver:1997:ASF} uniform expansions can be derived for all combinations of the parameters $n$, $\alpha$, $\beta$.

Let $\sigma$, $\tau$ and $\kappa$ be defined by
\begin{equation}\label{eq:Intro02}
\sigma=\frac{\alpha+\beta}{2\kappa},\quad \tau=\frac{\alpha-\beta}{2\kappa},\quad \kappa=n+\tfrac12(\alpha+\beta+1).
\end{equation}
Then $W(x)$ satisfies the differential equation
\begin{equation}\label{eq:Intro03}
\frac{d^2}{dx^2}W(x)=-\frac{\kappa^2(x_+-x)(x-x_-) +\frac14(x^2+3)}{(1-x^2)^2} W(x),
\end{equation}
where
\begin{equation}\label{eq:Intro04}
x_\pm=-\sigma\tau\pm\sqrt{(1-\sigma^2)(1-\tau^2)};
\end{equation}
$x_-$ and $x_+$ are called turning points. We have $-1\le x_-\le x_+\le 1$ when $\alpha$ and $\beta$ are positive.  When $\sigma^2+\delta^2=1$, one of the turning points $x_\pm$ is zero.

When we skip the term $\frac14(x^2+3)$ of the denominator 
in \eqref{eq:Intro03}, the differential equation becomes one for the Whittaker or Kummer functions, with special case the Laguerre polynomial, and when we take $\alpha=\beta$  the equation becomes a differential equation for the Gegenbauer polynomial.

When $\kappa$ is large we can make a few observations.

\begin{enumerate}
\item
If $n\gg \alpha+\beta$, then $\sigma\to0$ and $\tau\to0$. Hence, $x_-\to-1$ and $x_+\to1$. This is the standard case for large degree, the zeros are spread over the complete interval $(-1,1)$. 
\item
When $\alpha$ and/or $\beta$ become large as well, the zeros are inside the interval $(x_-,x_+)$. When, in addition, $\alpha/\beta\to0$, the zeros shift to the right, when $\beta/\alpha\to0$, they shift to the left. See also the limit in \eqref{eq:Intro09}.  The zeros become all positive when $x_-\ge0$. In that case $\sigma^2+\delta^2\ge1$.

\item
When $x$ is in a closed neighborhood around $x_-$ that does not contain $-1$ and $x_+$, an expansion in terms of Airy functions can be given. Similar for $x$ in a closed neighborhood around $x_+$ that does not contain $x_-$ and $1$. The points $x_\pm$ are called turning points of the equation in \eqref{eq:Intro03}.
\item
When $-1\le x\le x_-(1+a)<x_+$, with $a$ a fixed positive small number, an expansion in terms of Bessel functions can be given. Similar for $x_-<x_+(1-a)\le x\le 1$. The latter case corresponds to the limit
\begin{equation}\label{eq:Intro05}
\lim_{n\to\infty} n^{-\alpha}P_n^{(\alpha,\beta)}\left(1-\frac{x^2}{2n^2}\right)=\left(\frac{2}{x}\right)^\alpha J_\alpha(x).
\end{equation} 
Also, $\sqrt{x}J_\alpha\left(\alpha\sqrt{x}\right)$  satisfies the differential equation
\begin{equation}\label{eq:Intro06}
\frac{d^2}{dx^2}w(x)=\left(\alpha^2 \frac{1-x}{4x^2}-\frac{1}{4x^2}\right)w(x),
\end{equation}
in which $x=1$ is a turning point when $\alpha $ is large.

\item
If $ \alpha+\beta\gg n$, then $\sigma\to1$ and the turning points $x_-$ and $x_+$ coalesce at~$-\tau$. When $\alpha$ and $\beta$ are of the same order, the point $-\tau$ lies properly inside $(-1,1)$, and this case has been studied in \cite{Olver:1980:UAE} to obtain approximations of Whittaker functions in terms of parabolic cylinder functions. In the present case the parameters are such that the parabolic cylinder functions become Hermite polynomials.  This corresponds to the limit (see \cite{Lopez:1999:AOP})
\begin{equation}\label{eq:Intro07}
\lim_{\alpha,\beta\to\infty} \left(\frac{8}{\alpha+\beta}\right)^{n/2}\,
P_n^{(\alpha,\beta)}\left(x\sqrt{{\frac{2}{\alpha+\beta}}}-
\frac{\alpha-\beta}{\alpha+\beta}\right)=\frac1{n!}\,H_n(x),
\end{equation}
derived under the conditions 
\begin{equation}\label{eq:Intro08}
x=\bigO(1),\quad n=\bigO(1),\quad \frac{\alpha-\beta}{\alpha+\beta}=o(1),\quad \alpha, \beta\to\infty.
\end{equation}

\item
If $\alpha\gg\beta$, then $\tau\to1$, and $x_-$ and $x_+$ coalesce at $-\sigma$; if $\beta/\kappa=o(1)$, then the collision will happen at $-1$. Approximations in terms of Laguerre polynomials can be given. This corresponds to the limit
\begin{equation}\label{eq:Intro09}
\lim_{\alpha\to\infty}P^{(\alpha,\beta)}_{n}\bigl((2x/\alpha)-1\bigr)=(-1)^{n%
}L^{(\beta)}_{n}(x).
\end{equation}
Similar for $\beta\gg\alpha$, in which case $L^{(\alpha)}_{n}(x)$ becomes the approximant.

\end{enumerate}

As explained earlier, we consider in this paper the second case: new expansions of $P^{(\alpha,\beta)}_{n}(x)$, and its zeros and weights in terms of elementary functions. Preliminary results regarding the role of Gegenbauer and Laguerre polynomials as approximants  can be found in \cite{Temme:1990:PAE}.

\section{An integral representation and its saddle points}\label{sec:Jacnabelfunint}
The Rodrigues formula for the Jacobi polynomials reads  (see \cite[\S18.15(ii)]{Koornwinder:2010:OPS})
\begin{equation}\label{eq:int01}
P_n^{(\alpha,\beta)}(x)=\frac{(-1)^n}{2^n n!\,w(x)}\frac{d^n}{dx^n}\left(w(x)(1-x^2)^n\right),
\end{equation}
where
\begin{equation}\label{eq:int02}
w(x)=(1-x)^\alpha(1+x)^\beta.
\end{equation}
This gives the Cauchy integral representation
\begin{equation}\label{eq:int03}
P_n^{(\alpha,\beta)}(x)=\frac{(-1)^n}{2^n\,w(x)}\frac{1}{2\pi i}\int_\calC \frac{w(z)(1-z^2)^n}{(z-x)^{n+1}}\,dz, \quad x\in(-1,1),
\end{equation}
where the contour $\calC$ is a circle around the point $z=x$ with radius small enough to have the points $\pm1$ outside the circle. 

We write this in the form\footnote{The multi-valued functions of the integrand are discussed in Remark~\ref{rem:rem01}.}
\begin{equation}\label{eq:int04}
P_n^{(\alpha,\beta)}(x)=\frac{-1}{2^n\,w(x)}\frac{1}{2\pi i}\int_\calC e^{-\kappa \phi(z)}\,\frac{dz}{\sqrt{(1-z^2)(x-z)}}, 
\end{equation}
where 
\begin{equation}\label{eq:int05}
\kappa=n+\tfrac12(\alpha+\beta+1). 
\end{equation}
and 
\begin{equation}\label{eq:int06}
\phi(z)=-\frac{n+\alpha+\frac12}{\kappa}\ln(1-z)-\frac{n+\beta+\frac12}{\kappa}\ln(1+z)+\frac{n+\frac12}{\kappa}\ln(x-z).
\end{equation}
We introduce the notation
\begin{equation}\label{eq:int07}
\sigma=\frac{\alpha+\beta}{2\kappa},\quad \tau=\frac{\alpha-\beta}{2\kappa},
\end{equation}
and it follows that
\begin{equation}\label{eq:int08}
\phi(z)=-(1+\tau)\ln(1-z)-(1-\tau)\ln(1+z)+(1-\sigma)\ln(x-z).
\end{equation}
The saddle points $z_{\pm}$ follow from the zeros of
\begin{equation}\label{eq:int09}
\phi^\prime(z)= - \frac{(1+\sigma)z^2+2(\tau-x)z+1-\sigma-2\tau x}{(1-z^2)(x-z)},
\end{equation}
and are given by
\begin{equation}\label{eq:int10}
\begin{array}{@{}r@{\;}c@{\;}l@{}}
z_{\pm}&=&\dsp{\frac{x-\tau\pm iU(x)}{1+\sigma},}\\[8pt]
U(x)&=&\sqrt{1-2\sigma\tau x-\tau^2-\sigma^2-x^2}=\sqrt{(x_+-x)(x-x_-)},
\end{array}
\end{equation}
where (see also \eqref{eq:Intro04})
\begin{equation}\label{eq:int11}
x_{\pm}=-\sigma\tau\pm\sqrt{(1-\sigma^2)(1-\tau^2)}.
\end{equation}
In this representation we assume that $x_-\le x \le x_+$, in which $x$-domain the zeros  of the Jacobi polynomial are located.

\begin{remark}\label{rem:rem01}
The starting  integrand in \eqref{eq:int03} has a pole at $z=x$, while the one of \eqref{eq:int04} shows an algebraic singularity at $z=x$ and $\phi(z)$ defined in \eqref{eq:int06} has a logarithmic singularity at this point. To handle this from the viewpoint of multi-valued functions, we  can introduce a branch cut for the functions involved from $z=x$ to the left, assuming that the phase of $z-x$ is zero when $z>x$, equals $-\pi$ when $z$ approaches $-1$ on the lower part of the saddle point contour of the integral  in \eqref{eq:int04}, and $+\pi$ on the upper side. Because the saddle points $z_\pm$ stay off the interval $(-1,1)$, we do not need to consider function values on the branch cuts for the asymptotic analysis.
\eoremark
\end{remark}

\section{Deriving the asymptotic expansion}\label{sec:Jacnabelfun}
 We derive an expansion in terms of elementary functions which is valid for $x\in[x_-(1+\delta),x_+(1-\delta)]$, where $x_\pm$ are the turning points defined in \eqref{eq:int11} and $\delta$ is a fixed positive small number. Also, we assume that $\sigma\in[0,\sigma_0]$ and $\tau\in[-\tau_0,\tau_0]$, where $\sigma_0$ and $\tau_0$ are fixed positive numbers smaller than $1$. The case $\sigma\to1$ is explained in Case~5 of Section~\ref{sec:phen}. A similar phenomenon occurs when $\tau\to\pm1$.

 First we consider contributions from the saddle point $z_+$ using the  transformation
\begin{equation}\label{eq:Jacasymp01}
\phi(z)-\phi(z_+)=\tfrac12w^2
\end{equation}
for the contour from $z=+1$ to $z=-1$  through $z_+$, with $\phi(z)$ and $z_+$ given in \eqref{eq:int08} and 
\eqref{eq:int10}. This transforms the part of  the integral in \eqref{eq:int04} that runs with $\Im z\ge0$  into 
\begin{equation}\label{eq:Jacasymp02}
P^+=\frac{e^{-\kappa\phi(z_+)}}{2^n\,w(x)}\frac{1}{2\pi i}\int_{-\infty}^\infty e^{-\frac12\kappa w^2}f_+(w)\,dw,
\end{equation}
where
\begin{equation}\label{eq:Jacasymp03}
f_+(w)= \frac{1}{\sqrt{(1-z^2)(x-z)}}\frac{dz}{dw},\quad \frac{dz}{dw}=\frac{w}{\phi^\prime(z)}.
\end{equation}
We expand $\dsp{f_+(w)=\sum_{j=0}^\infty f_j^+w^j}$, where
\begin{equation}\label{eq:Jacasymp04}
f_0^+= \frac{1}{\sqrt{(1-z_+^2)(x-z_+)\phi^{\prime\prime}(z_+)}}=\Frac{e^{\frac14\pi i}}{\sqrt{2U(x)}},
\end{equation}
and $U(x)$ is defined in \eqref{eq:int10}. Because the contribution from the saddle point $z_-$ is the complex conjugate of that from $z_+$\footnote{We assume that $x\in(x_-,x_+)$ and that $\alpha$ and $\beta$ are positive.}, we take twice the real part of the contribution from $z_+$ and obtain the expansion
\begin{equation}\label{eq:Jacasymp05}
P_n^{(\alpha,\beta)}(x)\sim\Re\frac{e^{-\kappa\phi(z_+)-\frac14\pi i}}{2^{n}\,w(x)\sqrt{\pi \kappa U(x)}}\,\sum_{j=0}^\infty \frac{c_{j}^+}{\kappa^j}, \quad c_{j}=2^j\left(\tfrac12\right)_j \frac{f_{2j}^+}{f_0^+}.
\end{equation}

Evaluating $\phi(z_+)$ we find
\begin{equation}\label{eq:Jacasymp06}
\begin{array}{@{}r@{\;}c@{\;}l@{}}
\phi(z_+)&=&-\ln 2+\psi+\xi+i\chi(x),\\[8pt]
\psi&=&-\frac12(1-\tau)\ln(1-\tau)-\frac12(1+\tau)\ln(1+\tau)\ +\\[8pt]
&&\frac12(1+\sigma)\ln(1+\sigma)+\frac12(1-\sigma)\ln(1-\sigma),\\[8pt]
\xi(x)&=&-\frac12(\sigma+\tau)\ln(1-x)-\frac12(\sigma-\tau)\ln(1+x),\\[8pt]
\chi(x)&=&\dsp{(\tau+1)\arctan\frac{U(x)}{1-x+\sigma+\tau}+(\tau-1)\arctan\frac{U(x)}{1+x+\sigma-\tau}\ +}\\[8pt]
&&(1-\sigma)\,{\rm{atan}}2(-U(x),\tau+x\sigma).
\end{array}
\end{equation}

\begin{figure}
\begin{center}
\epsfxsize=5cm \epsfbox{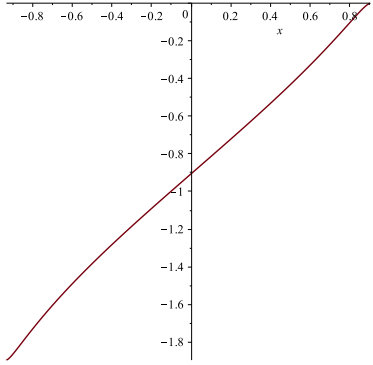}
\caption{
\label{fig:fig01} The quantity $\chi(x)$ defined in \eqref{eq:Jacasymp06} for $x\in(x_-,x_+)$; $\alpha=90$, $\beta=75$, $n=125$. For these values, $\kappa=208$, $\sigma=\frac{165}{416}$, $\tau=\frac{15}{416}$, $x_-=-0.931$, $x_+=0.903$.}
\end{center}
\end{figure}

In Figure~\ref{fig:fig01} we show a graph of $\chi(x)$ on $(x_-,x_+)$ for $\alpha=90$, $\beta=75$, $n=125$. For these values, $\kappa=208$, $\sigma=\frac{165}{416}$, $\tau=\frac{15}{416}$, $x_-=-0.931$, $x_+=0.903$. At the left endpoint we have $\chi(x_-)=-(1-\sigma)\pi=-1.896$.

\begin{remark}\label{rem:rem02}
The denominators of the first and second arctan functions of $\chi(x)$  in \eqref{eq:Jacasymp06} are always positive on  $(x_-,x_+)$; this follows easily from the relations in \eqref{eq:int07}. The function ${\rm{atan}}2(y,x)$ in the third term of $\chi(x)$ denotes the phase $\in(-\pi,\pi]$ of the complex number $x+iy$. Because $\tau+x\sigma$ may be negative  on $(x_-,x_+)$ we cannot use the standard arctan function for that term.
\eoremark
\end{remark}

Observe that $e^{-\kappa\xi(x)}=\sqrt{w(x)}$, with $w(x)$ defined in \eqref{eq:int02}. To compute $x$ from $\chi(x)$, for example by using a Newton-procedure, it is convenient to know that
\begin{equation}\label{eq:Jacasymp07}
\frac{d\chi(x)}{dx}=\frac{U(x)}
{\left(1-x^2\right)}.
\end{equation}

We return to the result in \eqref{eq:Jacasymp05} and split the coefficients of \eqref{eq:Jacasymp05} in real and imaginary parts. We write $c_j^+=p_j+iq_j$, and obtain
\begin{equation}\label{eq:Jacasymp08}
\begin{array}{@{}r@{\;}c@{\;}l@{}}
P_n^{(\alpha,\beta)}(x)&=&\dsp{\frac{2^{\frac12(\alpha+\beta+1)}e^{-\kappa\psi}}
{\sqrt{\pi \kappa w(x)U(x)}}W(x)},\\[8pt]
W(x)&=&\dsp{\cos\left(\kappa\chi(x)+\tfrac14\pi\right)P(x)+\sin\left(\kappa\chi(x)+\tfrac14\pi\right)Q(x),}
\end{array}
\end{equation}
with expansions
\begin{equation}\label{eq:Jacasymp09}
P(x)\sim \sum_{j=0}^\infty \frac{p_{j}}{\kappa^j},\quad Q(x)\sim \sum_{j=0}^\infty \frac{q_{j}}{\kappa^j}.
\end{equation}
Because $c_0^+=1$, we have $p_0=1$, $q_0=0$.

To evaluate the coefficients $f_{2j}^+$ of the expansion in \eqref{eq:Jacasymp05}, we need the coefficients $z_j^+$ of the expansion $z=z_++\sum_{j=1}^\infty z_j^+ w^j$ that follow from \eqref{eq:Jacasymp01}. The first values are 
\begin{equation}\label{eq:Jacasymp10}
\begin{array}{@{}r@{\;}c@{\;}l@{}}
 z_2^+&=&-\tfrac16 z_1^4\phi_3,\quad z_3^+=\dsp{\tfrac1{72}z_1^5\left(5z_1^2\phi_3^2-3\phi_4\right)},\\[8pt]
z_4^+&=&\dsp{-\tfrac1{1080}z_1^6\left(9\phi_5-45z_1^2\phi_3\phi_4+40z_1^4\phi_3^2\right),}
\end{array}
\end{equation}
where $z_1=z_1^+=1/\sqrt{\phi^{\prime\prime}(z_+)}$ and $\phi_j$ denotes the $j$th derivative of $\phi(z)$ at the saddle point $z=z_+$ defined in \eqref{eq:int10}.

With these coefficients we expand $f(w)$ defined in \eqref{eq:Jacasymp04}. This gives
\begin{equation}\label{eq:Jacasymp11}
\begin{array}{@{}r@{\;}c@{\;}l@{}}
c_1^+&=&\dsp{-\frac{ z_+}{8z_1 (1- z_+^2)^2 (x- z_+)^2}}\Bigl(-6 z_1^3  z_+^2+3 z_1^3-72 z_1 z_2  z_+^2 x\ +\\[8pt]
&&24 z_1  z_+ z_2 x^2-24 z_1  z_+^3 z_2 x^2-48 z_3 x  z_+-48 z_3  z_+^2 x^2+96 z_3  z_+^3 x\ +\\[8pt]
&&24 z_3  z_+^4 x^2-48 z_3  z_+^5 x-12 z_1  z_+ z_2+48 z_1  z_+^3 z_2-48 z_3  z_+^4\ +\\[8pt]
&&24 z_3  z_+^6+12 z_1 z_2 x-36 z_1 z_2  z_+^5-4 z_1^3 x  z_++8 z_1^3  z_+^2 x^2-20 z_1^3  z_+^3 x\ + \\[8pt]
&&4 z_1^3 x^2+15 z_1^3  z_+^4+24 z_3 x^2+24 z_3  z_+^2+60 z_1 z_2  z_+^4 x\Bigr),
\end{array}
\end{equation}
where $z_j$ denotes $z_j^+$. The coefficients $p_1$ and $q_1$ of the expansions in \eqref{eq:Jacasymp09} follow from $c_1^+=p_1+iq_1$.

\subsection{Expansion of the derivative}\label{sec:Pderiv}
For the weights of the Gauss quadrature it is convenient to have an expansion of $\dsp{\frac{d}{dx}}P_n^{(\alpha,\beta)}(x)$. Of course this follows from using \eqref{eq:Jacasymp08} with different values of $\alpha$ and $\beta$ and the relation
\begin{equation}\label{eq:Jacasymp12}
\frac{d}{dx}P_n^{(\alpha,\beta)}(x)=\tfrac{1}{2}\left(\alpha+\beta+n+1\right)P_{n-1}^{(\alpha+1,\beta+1)}(x),
\end{equation}
but it is useful to have a representation in terms of the same parameters.

By straightforward differentiation of \eqref{eq:Jacasymp08} we obtain 
\begin{equation}\label{eq:Jacasymp13}
\begin{array}{@{}r@{\;}c@{\;}l@{}}
\dsp{\frac{d}{dx}P_n^{(\alpha,\beta)}(x)}&=&
\dsp{-\sqrt{\frac{\kappa}{\pi}}\,2^{\frac12(\alpha+\beta+1)e^{-\kappa\psi}}\chi^\prime(x)A(x)
 \ \times}\\[8pt]
&&\dsp{\left(\sin\left(\kappa\chi(x)+\tfrac14\pi\right)R(x)-\cos\left(\kappa\chi(x)+\tfrac14\pi\right)S(x)\right)},
\end{array}
\end{equation}
where $\chi^\prime(x)$ is given in \eqref{eq:Jacasymp07} and
\begin{equation}\label{eq:Jacasymp14}
\begin{array}{@{}r@{\;}c@{\;}l@{}}
A(x)&=&\dsp{\frac{1}{\sqrt{w(x)U(x)}}},\\[8pt]
R(x)&=&\dsp{P(x)-\frac{1}{\kappa\chi^\prime(x)}Q^\prime(x)-\frac{A^\prime(x)}{\kappa A(x)\chi^\prime(x)}Q(x)},\\[8pt]
S(x)&=&\dsp{Q(x)+\frac{1}{\kappa\chi^\prime(x)}P^\prime(x)+\frac{A^\prime(x)}{\kappa A(x)\chi^\prime(x)}P(x)}.
\end{array}
\end{equation}
We have the expansions
\begin{equation}\label{eq:Jacasymp15}
R(x)\sim \sum_{j=0}^\infty \frac{r_{j}}{\kappa^j},\quad S(x)\sim \sum_{j=0}^\infty \frac{s_{j}}{\kappa^j},
\end{equation}
where the coefficients follow  from the relations in \eqref{eq:Jacasymp14}. The first coefficients are $r_0=p_0=1$, $s_0=q_0=0$, and 
\begin{equation}\label{eq:Jacasymp16}
r_1=p_1,\quad s_1=q_1+\frac{A^\prime(x)}{A(x)\chi^\prime(x)}.
\end{equation}

\section{Expansion of the zeros}\label{sec:Jacnabzer}
A zero $x_\ell$, $1\le \ell\le n$, of $P_n^{(\alpha,\beta)}(x)$ follows from the zeros of (see \eqref{eq:Jacasymp08})
\begin{equation}\label{eq:Jaczeros01}
W(x)=\cos\left(\kappa\chi(x)+\tfrac14\pi\right)P(x)+\sin\left(\kappa\chi(x)+\tfrac14\pi\right)Q(x),
\end{equation}
where $\chi(x)$ is defined in \eqref{eq:Jacasymp08}. For a first approximation we put the cosine term equal to zero. That is, we can write
\begin{equation}\label{eq:Jaczeros02}
\kappa\chi(x)+\tfrac14\pi=\tfrac12\pi-(n+1-\ell)\pi,
\end{equation}
where $\ell$ is some integer. It appears that this choice in the right-hand side is convenient for finding the $\ell$th zero.

Because the expansions in \eqref{eq:Jacasymp09} are valid for $x$ properly inside $(x_-,x_+)$, we may expect that the approximations of the zeros in the middle of this interval will be much better than those near the endpoints. We describe how to compute approximations of all  $n$ zeros by considering the zeros of $\cos\left.(\chi(x)\kappa +\frac14\pi\right)$.

We start with $\ell=1$ and using \eqref{eq:Jaczeros02} we compute $\chi_1=\left(\frac14-n\right)\pi/\kappa$. Next we compute an approximation of the zero $x_1$ by inverting the equation $\chi(x)=\chi_1$, where $\chi(x)$ is defined in \eqref{eq:Jacasymp08}. For a Newton procedure we can use $x_-+1/n$ as a starting value.

\begin{example}\label{exemp:ex01}
When we take $\alpha=50$, $\beta=41$, $n=25$, we have $\kappa=71$, $\sigma=91/142$, $\tau=9/142$. We find $\chi_1= -1.095133$ and the starting value of the Newton procedure is $x= -0.7667437$. We find  $x_1\doteq -0.7415548$. Comparing this with the first zero computed by using the solver of Maple to compute the zeros of the Jacobi polynomial with Digits = 16, we find a relative error $0.00074$. 

For the next zero $x_2$, we compute $\chi_2$ from \eqref{eq:Jaczeros02} with $\ell=2$, use $x_1$  as a starting value for the Newton procedure, and find $x_2 \doteq-0.682106$, with relative error $0.00032$. And so on.
The best result is for $x_{13}$ with relative error $0.000013$, and the worst result is for $x_{25}$ with a relative error  
$0.0010$.
\eoexample
\end{example}
\begin{remark}\label{rem:rem03}
We don't have a proof that the found zero always corresponds with the $\ell$th zero, when we start with \eqref{eq:Jaczeros02}. In a number of tests we have found all agreement with this choice. 
\eoremark
\end{remark}

To obtain higher approximations of the zeros, we use the method described in our earlier papers.  We assume that the zero $x_\ell$ has an asymptotic expansion
\begin{equation}\label{eq:Jaczeros03}
x_\ell=\xi_0+\eps,\quad \eps\sim  \frac{\xi_2}{\kappa^2}+\frac{\xi_4}{\kappa^4}+\ldots,
\end{equation}
where $\xi_0$ is the value obtained as a first approximation by the method just described. 

The function $W(x)$ defined in \eqref{eq:Jaczeros01} can be expanded  at $\xi_0$ and we have
\begin{equation}\label{eq:Jaczeros04}
W(x_\ell)=W(\xi_0+\eps)=W(\xi_0)+\frac{\eps}{1!}W^\prime(\xi_0)+ \frac{\eps^2}{2!}W^{\prime\prime}(\xi_0)+\ldots = 0,
\end{equation}
where the derivatives are with respect to $x$.  We find upon substituting the expansions of $\eps$ and those of $P$ and $Q$ given  \eqref{eq:Jacasymp09},  and comparing equal powers of $\kappa$, that the first coefficients are
\begin{equation}\label{eq:Jaczeros05}
\begin{array}{@{}r@{\;}c@{\;}l@{}}
\xi_2&=&\dsp{\frac{\left(1-x^2\right)q_1(x)}{U(x)}},\\[8pt]
\xi_4&=&\dsp{\frac{1}{6U(x)^4}}\Bigl(3x^5q_1^2+3x^4q_1^2\sigma\tau-6x^3q_1^2-6x^2q_1^2\sigma\tau+3q_1^2x+3q_1^2\sigma\tau\ +\\[8pt]
&&\bigl(6q_1^\prime q_1x^4+6x^3q_1^2-12q_1^\prime x^2q_1-6xq_1^2+6q_1^\prime q_1\bigr)U(x)^2\ +\\[8pt]
&&\bigl(6p_2x^2q_1+2q_1^3x^2+6q_3-6p_2q_1-6q_3x^2-2q_1^3\bigr)U(x)^3\Bigr),
\end{array}
\end{equation}
where $U(x)$ is defined in \eqref{eq:int10}, and $x$ takes the value of the first approximation of the zero as obtained in Example~\ref{exemp:ex01}.

When we take the same values  $\alpha=50$, $\beta=41$, $n=25$ as in Example~\ref{exemp:ex01}, and use \eqref{eq:Jaczeros03} with the term $\xi_2/\kappa^2$  included, we obtain for the zero $x_{13}$  a relative error $0.80\times10^{-9}$. With also  the term $\xi_4/\kappa^4$  included we find for $x_{13}$ a relative error $0.13\times10^{-12}$.

\begin{figure}
\epsfxsize=15cm \epsfbox{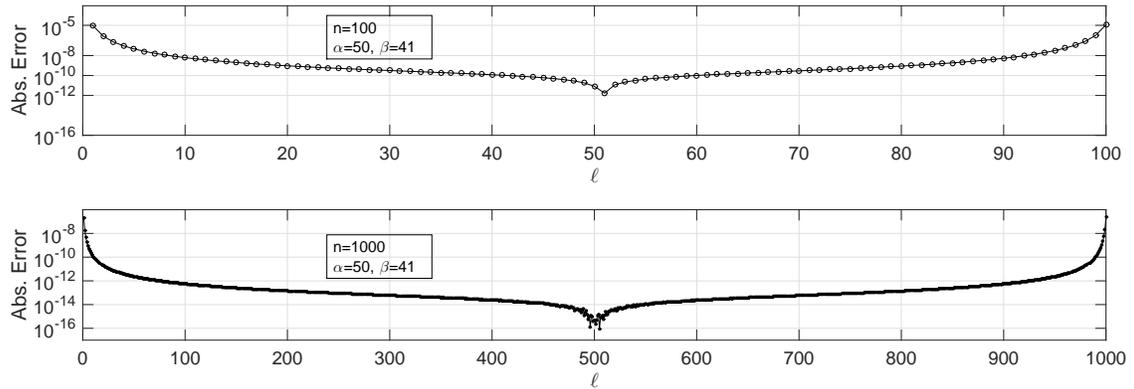}
\caption{
\label{fig:fig02} Performance of the asymptotic expansion for computing the zeros of $P^{(\alpha,\beta)}_n(x)$ for
 $\alpha=50$, $\beta=41$ and $n=100,\,1000$.}
\end{figure}

A more extensive test of the expansion is shown in Figure ~\ref{fig:fig02}.
The label $\ell$ in the abscissa represents the order of
the zero (starting from $\ell = 1$ for the smallest zero).
In this figure we compare the approximations to
the zeros obtained with the asymptotic expansion against the results of a Maple
implementation (with a large number of digits) of an iterative algorithm which uses
the global fixed point method of \cite{Segura:2010:RCO}. The Jacobi polynomials used in this algorithm
are computed by using the intrinsic Maple function. As before, we use \eqref{eq:Jaczeros03} with the term $\xi_2/\kappa^2$  included.
As can be seen, for $n = 100$ the use of the expansion allows the computation
of  the zeros $x_\ell$, $10\le\ell\le90$, with absolute error less than $10^{-8}$. When $n = 1000$, 
an absolute accuracy better than $10^{-12}$ can be obtained for about 90\% of the zeros of
the Jacobi polynomials. The results become less accurate for the zeros near the endpoints $\pm1$, as expected.

\begin{figure}
\epsfxsize=15cm \epsfbox{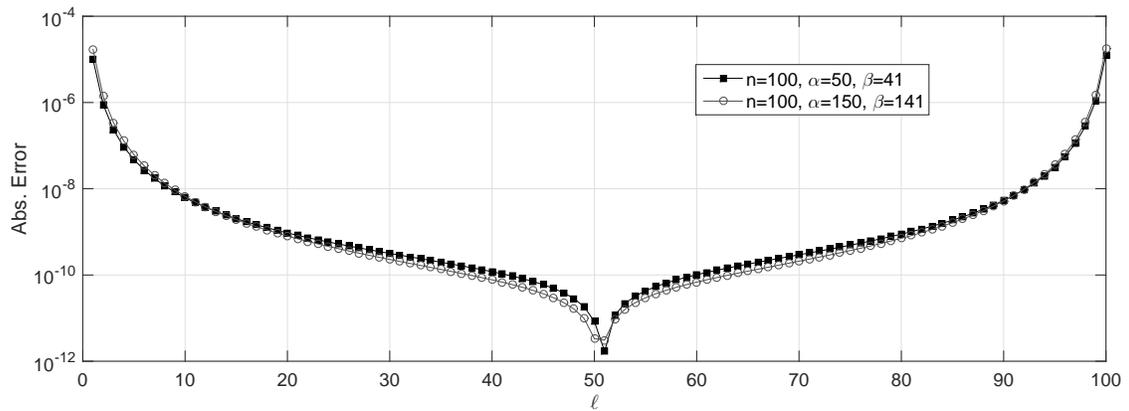}
\caption{
\label{fig:fig03} Performance of the asymptotic expansion for computing the zeros for $n=100$ and $\alpha=50$, $\beta=41$ compared with $\alpha=150$, $\beta=141$. }
\end{figure}

In Figure~\ref{fig:fig03} we show the absolute errors for $n=100$ and $\alpha=50$, $\beta=41$ compared with $\alpha=150$, $\beta=141$. We see that the accuracy is slightly better for the larger parameters, and that the asymptotics is quite uniform when $\alpha$ and $\beta$ assume larger values.

\section{The weights of the Gauss-Jacobi quadrature}\label{sec:weights}

As we did in \cite{Gil:2019:NIG}, and in our earlier paper \cite{Gil:2018:GHL} for the Gauss--Hermite and Gauss--Laguerre quadratures, it is convenient 
to introduce scaled weights. In terms of the derivatives of the Jacobi polynomials, the classical form of the weights of the Gauss-Jacobi quadrature can be written as
\begin{equation}\label{eq:weights01}
\begin{array}{@{}r@{\;}c@{\;}l@{}}
w_\ell  &=&
\dsp{\frac{M_{n,\alpha,\beta}}{ \left(1-x_\ell^2\right)  P_n ^{(\alpha ,\beta)\prime}(x_\ell)^2}},
 \\
&&\\
M_{n,\alpha,\beta}&=&\dsp{2^{\alpha+\beta+1}\frac{\Gamma (n+\alpha+1)\Gamma (n+\beta+1)}{n! \Gamma (n+\alpha+\beta+1 )}}.
\end{array}
\end{equation}

In Figure 4 we show the relative errors in the computation of the weights
$w_\ell$ defined in \eqref{eq:weights01}, with the derivative of the Jacobi polynomial computed by
using the relation in \eqref{eq:Jacasymp12}. We have used the representation in \eqref{eq:Jacasymp08}, with the
asymptotic series \eqref{eq:Jacasymp09} truncated after $j = 3$ and the expansion  \eqref{eq:Jaczeros03} for the
nodes   with the term $\xi_2/\kappa^2$ included. The relative errors are obtained
by using high-precision results computed by using Maple.

\begin{figure}
\epsfxsize=15cm \epsfbox{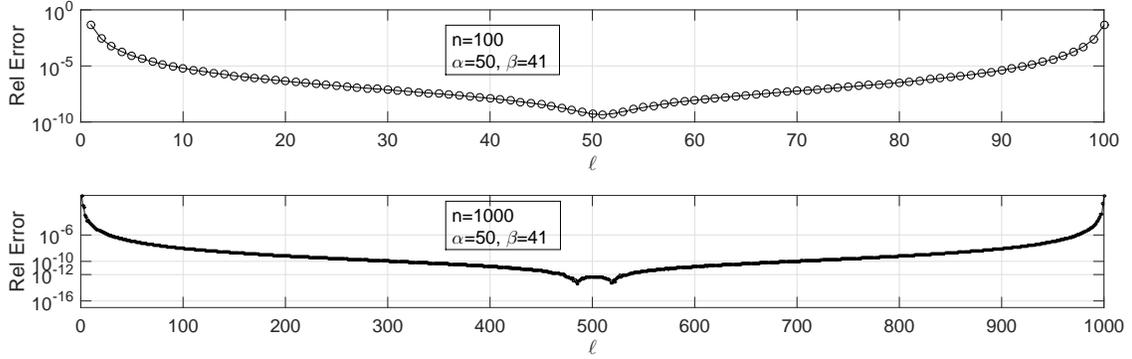}
\caption{
\label{fig:fig04} Performance of the computation of the weights $w_\ell$ by using the asymptotic expansion of the Jacobi polynomial for
 $\alpha=50$, $\beta=41$ and $n=100,\,1000$.}
\end{figure}

As an alternative we consider the scaled weights defined by
\begin{equation}\label{eq:weights02}
\omega_\ell=\frac{1}{v^\prime(x_\ell)^{2}},
\end{equation}
where
\begin{equation}\label{eq:weights03}
v(x)=C_{n,\alpha,\beta}\,
(1-x)^{a}(1+x)^{b} P_n^{(\alpha,\beta)}(x),
\end{equation}
and we choose $a$ and $b$ such that $v^{\prime\prime}(x_\ell)=0$; $C_{n,\alpha,\beta}$ does not depend on $x$, and will be chosen later. 
We have
\begin{equation}\label{eq:weights04}
\begin{array}{@{}r@{\;}c@{\;}l@{}}
v^\prime(x)&=&C_{n,\alpha,\beta}\bigl(
\left(-a(1-x)^{a-1}(1+x)^{b} +b(1-x)^{a}(1+x)^{b-1}\right)P_n^{(\alpha,\beta)}(x)\ +\\[8pt]
&&(1-x)^{a}(1+x)^{b}P_n^{(\alpha,\beta)\prime}(x)\bigr).
\end{array}
\end{equation}
Evaluating $v^{\prime\prime}(x_\ell)$, we find 
\begin{equation}\label{eq:weights05}
\begin{array}{@{}r@{\;}c@{\;}l@{}}
v^{\prime\prime}(x_\ell)&=&C_{n,\alpha,\beta}(1-x_\ell)^{a}(1+x_\ell)^{b}(1-x_\ell^2)\ \times\\[8pt]
&&\left((1-x_\ell^2)P_n^{(\alpha,\beta)\prime\prime}(x_\ell)+
2\left(b-a-(a+b)x_\ell\right)P_n^{(\alpha,\beta)\prime}(x_\ell)\right),
\end{array}
\end{equation}
where we skip the term containing $P_n^{(\alpha,\beta)}(x_\ell)$, because $x_\ell$ is a zero.

The differential equation of the Jacobi polynomials is
\begin{equation}\label{eq:weights06}
\left(1-x^2\right) y^{\prime\prime}(x)+\left(\beta-\alpha-(\alpha+\beta+2)x\right)y^\prime(x)+n(\alpha+\beta+n+1)y(x)=0,
\end{equation}
and we see that $v^{\prime\prime}(x_\ell)=0$ if we take $a=\frac12(\alpha+1)$, $b=\frac12(\beta+1)$. 

We obtain
\begin{equation}\label{eq:weights07}
v(x)=C_{n,\alpha,\beta}\,
(1-x)^{\frac12(\alpha+1)}(1+x)^{\frac12(\beta+1)} P_n^{(\alpha,\beta)}(x),
\end{equation}
with properties
\begin{equation}\label{eq:weights08}
v^\prime(x_\ell)=C_{n,\alpha,\beta}\,
(1-x_\ell)^{\frac12(\alpha+1)}(1+x_\ell)^{\frac12(\beta+1)} P_n^{(\alpha,\beta)\prime}(x_\ell),
\quad v^{\prime\prime}(x_\ell)=0.
\end{equation}
The weights $w_\ell$ are related with the scaled weights $\omega_\ell$ by
\begin{equation}\label{eq:weights09}
w_\ell =M_{n,\alpha,\beta}C^2_{n,\alpha,\beta}(1-x_\ell)^{\alpha}(1+x_\ell)^{\beta} \omega_\ell.
\end{equation}

The advantage of computing scaled weights is that, similarly as described in \cite{Gil:2018:GHL},
scaled weights do not underflow/overflow for large parameters. In additional, they are well-conditioned
as a function of the roots $x_\ell$. Indeed, introducing the notation
\begin{equation}\label{eq:weights10}
V(x)=\frac{1}{v^\prime(x)^{2}},
\end{equation}
the scaled weights are $\omega_\ell=V(x_\ell)$ and $V^\prime(x_\ell)=0$ because $v^{\prime\prime}(x_\ell)=0$. The vanishing derivative of $V(x)$ at $x_\ell$ may result in a more accurate numerical evaluation of the scaled weights.

When considering the representation of the Jacobi polynomials in \eqref{eq:Jacasymp08},
the function $v(x)$ can be written as 
\begin{equation}\label{eq:weights11}
v(x)=
\frac{2^{\frac12(\alpha+\beta+1)}}{\sqrt{\pi \kappa}}\,C_{n,\alpha,\beta}e^{-\kappa\psi}Z(x)W(x),\quad Z(x)=\sqrt{\frac{1-x^2}{U(x)}},
\end{equation}
where $U(x)$ is defined in \eqref{eq:int10}. For scaling $v(x)$ we choose
\begin{equation}\label{eq:weights12}
C_{n,\alpha,\beta}=2^{-\frac12(\alpha+\beta+1)}e^{\kappa\psi}.
\end{equation}
This gives
\begin{equation}\label{eq:weights13}
v(x)=
\frac{Z(x)W(x)}{\sqrt{\pi \kappa}}.
\end{equation}

For the numerical computation of $\psi$ defined in  \eqref{eq:Jacasymp06} for small values of $\sigma$ or $\tau$, we can use the expansion
\begin{equation}\label{eq:weights14}
(1-x)\ln(1-x)+(1+x)\ln(1+x)=\sum_{k=1}^\infty\frac{x^{2k}}{k(2k-1)},\quad \vert x\vert < 1.
\end{equation}

For computing the modified  Gauss weights it is convenient to have an expansion of the derivative of the function $v(x)$ of \eqref{eq:weights13},  with $W(x)$ defined in \eqref{eq:Jacasymp08} and $Z(x)$ in \eqref{eq:weights11}.

We have
\begin{equation}\label{eq:weights20}
\frac{d}{dx}v(x)=-\sqrt{\frac{\kappa}{\pi}}\chi^\prime(x)Z(x)
\left(\sin\left(\kappa\chi(x)+\tfrac14\pi\right)M(x)-\cos\left(\kappa\chi(x)+\tfrac14\pi\right)N(x)\right),
\end{equation}
where $\chi^\prime(x)$ is given in \eqref{eq:Jacasymp07} and
\begin{equation}\label{eq:weights21}
\begin{array}{@{}r@{\;}c@{\;}l@{}}
M(x)&=&\dsp{P(x)-\frac{1}{\kappa}p(x)Q^\prime(x)-\frac{1}{\kappa}q(x)Q(x)},\\[8pt]
N(x)&=&\dsp{Q(x)+\frac{1}{\kappa} p(x)P^\prime(x)+\frac{1}{\kappa}q(x)P(x)},
\end{array}
\end{equation}
where
\begin{equation}\label{eq:weights22}
\begin{array}{@{}r@{\;}c@{\;}l@{}}
p(x)&=&\dsp{\frac{1}{\chi^\prime(x)}=\frac{1-x^2}{U(x)},}\\[8pt]
q(x)&=&\dsp{\frac{Z^\prime(x)}{Z(x)\chi^\prime(x)}=\frac{(1-x^2)(x+\sigma\tau)-2xU^2(x)}{2U^3(x)}}.
\end{array}
\end{equation}
We have the expansions
\begin{equation}\label{eq:weights23}
M(x)\sim \sum_{j=0}^\infty \frac{m_{j}}{\kappa^j},\quad N(x)\sim \sum_{j=0}^\infty \frac{n_{j}}{\kappa^j},
\end{equation}
where the coefficients follow  from the relations in \eqref{eq:Jacasymp14}. The first coefficients are $m_0=p_0=1$, $n_0=q_0=0$, and for $j=1,2,3,\ldots$
\begin{equation}\label{eq:weights24}
\begin{array}{@{}r@{\;}c@{\;}l@{}}
m_j&=&\dsp{p_j-p(x)q_{j-1}^\prime-q(x)q_{j-1},}\\[8pt]
n_j&=&\dsp{q_j+p(x)p_{j-1}^\prime+q(x)p_{j-1}.}
\end{array}
\end{equation}

\begin{figure}
\epsfxsize=15cm \epsfbox{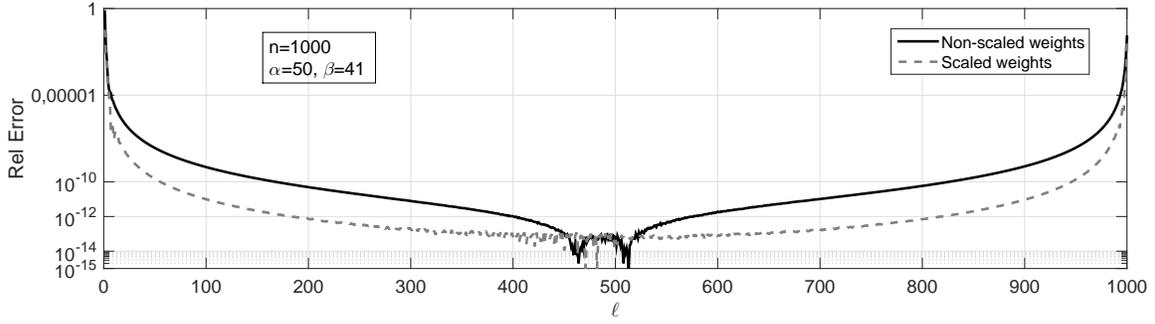}
\caption{
\label{fig:fig05} Comparison of the performance of the asymptotic expansions 
for computing non-scaled \eqref{eq:weights01} and scaled \eqref{eq:weights02} weights for
 $\alpha=50$, $\beta=41$ and $n=1000$.}
\end{figure}

As an example, Figure~\ref{fig:fig05} shows the performance of the asymptotic expansion \eqref{eq:weights20}
for computing the scaled weights  \eqref{eq:weights02} for  $\alpha=50$, $\beta=41$ and $n=1000$.
The computation of the non-scaled weights \eqref{eq:weights01} is shown as comparison.

\begin{figure}
\begin{center}
\epsfxsize=15cm \epsfbox{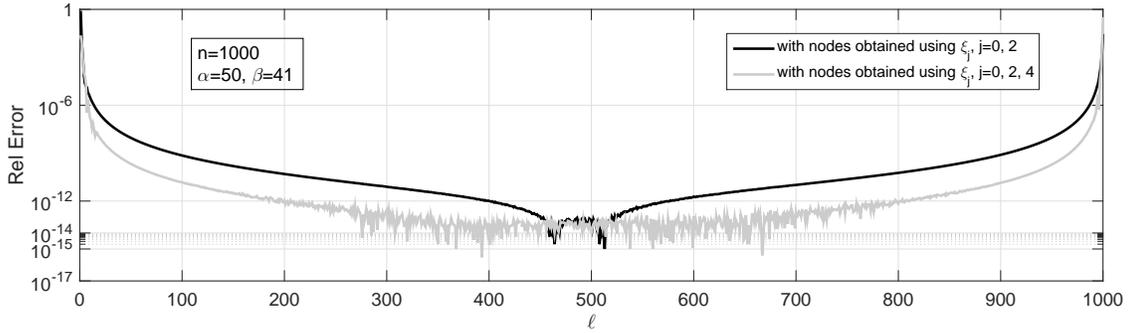}
\caption{
\label{fig:fig06} Performance of the computation of the weights $w_\ell$  defined in \eqref{eq:weights01} by using the asymptotic expansion of the Jacobi polynomial for
 $\alpha=50$, $\beta=41$ and $n=1000$. The comparison is between the expansion of the zeros in \eqref{eq:Jaczeros03} with the term  $\xi_4/\kappa^4$ included or not included.}
\end{center}
\end{figure}

\begin{figure}
\epsfxsize=15cm \epsfbox{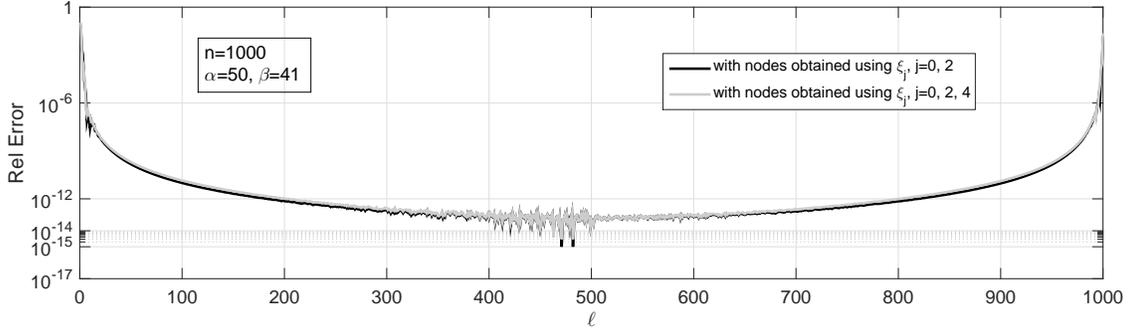}
\caption{
\label{fig:fig07} Same as in Figure~\ref{fig:fig06} for the scaled  weights $\omega_\ell$  defined in \eqref{eq:weights02}.}
\end{figure}

In Figure~\ref{fig:fig06} and  Figure~\ref{fig:fig07} we compare the effect of computing the weights  $w_\ell$ defined in \eqref{eq:weights01} and  the scaled weights $\omega_\ell$  defined in \eqref{eq:weights02} when we compute these weights with the asymptotic expansion of the zeros in \eqref{eq:Jaczeros03} with the term  $\xi_4/\kappa^4$ included or not included. From these computations it follows that the that the scaled weights are well-conditioned as a function of the nodes and therefore they are not so critically dependent on the accuracy of the nodes. Contrary the non-scaled weights have worse condition and the accuracy of the nodes is more important.

\subsection{About quantities appearing in the weights. }\label{sec:weightscoeff}

First we consider the term $e^{\kappa\psi}$, with $\psi$ given in \eqref{eq:Jacasymp06}. Using the relations in \eqref{eq:int07}, we have
\begin{equation}\label{eq:weights15}
\begin{array}{lll}
&\kappa(1+\tau)=n+\alpha+\frac12,&\kappa(1-\tau)=n+\beta+\frac12,\\[8pt]
&\kappa(1+\sigma)=n+\alpha+\beta+\frac12,&\kappa(1-\sigma)=n+\frac12,
\end{array}
\end{equation}
and this gives
\begin{equation}\label{eq:weights16}
\begin{array}{@{}r@{\;}c@{\;}l@{}}
e^{2\kappa\psi}&=&\dsp{\frac{\left(n+\alpha+\beta+\frac12\right)^{n+\alpha+\beta+\frac12}
\left(n+\frac12\right)^{n+\frac12}}
{\left(n+\alpha+\frac12\right)^{n+\alpha+\frac12}
\left(n+\beta+\frac12\right)^{n+\beta+\frac12}}}\\[8pt]
&=&\dsp{\frac{\Gamma\left( n+\alpha+\beta+\frac12\right)\Gamma\left(n
+\frac12 \right)}{\Gamma\left(n+\alpha+\frac12 \right)\Gamma\left( n+\beta+\frac12\right)}
\ \frac
{\Gamma^*\left(n+\alpha+\frac12 \right)\Gamma^*\left( n+\beta+\frac12\right)}
{\Gamma^*\left( n+\alpha+\beta+\frac12\right)\Gamma^*\left(n+\frac12 \right)}
 \times}\\[8pt]
&&\dsp{\sqrt{\frac{\left( n+\alpha+\beta+\frac12\right)\left( n+\frac12\right)}{\left( n+\alpha+\frac12\right)\left( n+\beta+\frac12\right)}}},
\end{array}
\end{equation}
where
\begin{equation}\label{eq:weights17}
\Gamma^*(z)=\sqrt{{z/(2\pi)}}\,e^z
z^{-z}\Gamma(z),\quad \phase\,z\in(-\pi,\pi),\quad z\ne0.
\end{equation}
We have $\Gamma^*(z)=1+\bigO(1/z)$ as $z\to\infty$.

It follows that (see \eqref{eq:weights01} and \eqref{eq:weights09})
\begin{equation}\label{eq:weights18}
\begin{array}{@{}r@{\;}c@{\;}l@{}}
M_{n,\alpha,\beta}C^2_{n,\alpha,\beta}&=&
\dsp{ \frac
{\Gamma\left(n+\alpha+1 \right)\Gamma\left( n+\beta+1\right)\Gamma\left( n+\alpha+\beta+\frac12\right)\Gamma\left(n+\frac12 \right)}
{\Gamma\left(n+\alpha+\frac12 \right)\Gamma\left( n+\beta+\frac12\right)\Gamma\left( n+\alpha+\beta+1\right)\Gamma\left(n+1 \right)}}\  \times\\[8pt]
&&\dsp{ 
\frac
{\Gamma^*\left(n+\alpha+\frac12 \right)\Gamma^*\left( n+\beta+\frac12\right)}
{\Gamma^*\left( n+\alpha+\beta+\frac12\right)\Gamma^*\left(n+\frac12 \right)}
\sqrt{\frac{\left( n+\alpha+\beta+\frac12\right)\left( n+\frac12\right)}{\left( n+\alpha+\frac12\right)\left( n+\beta+\frac12\right)}}.
}\end{array}
\end{equation}
Using $\Gamma\left(z+\frac12\right)/\Gamma(z)\sim z^{\frac12}$ as $z\to\infty$, we see that, in the case that $\alpha$, $\beta$ and $n$ are all large, we have $M_{n,\alpha,\beta}C^2_{n,\alpha,\beta}\sim1$, and that, when using more details on expansions of gamma functions and ratios thereof (see \cite[\S6.5]{Temme:2015:AMI}), we can obtain 
\begin{equation}\label{eq:weights19}
M_{n,\alpha,\beta}C^2_{n,\alpha,\beta}\sim 1+\frac{\sigma^2-\tau^2}{12(1-\sigma^2)(1-\tau^2)\kappa}+
\frac{(\sigma^2-\tau^2)^2}{288(1-\sigma^2)^2(1-\tau^2)^2\kappa^2}+\ldots,
\end{equation}
again, when $\alpha$, $\beta$ and $n$ are all large.

As observed in the first lines of Section~\ref{sec:Jacnabelfun}, in the present asymptotics we assume that $\sigma$ and $\vert\tau\vert$ are bounded away from $1$.

\section*{Acknowledgments}

We acknowledge financial support from Ministerio de Ciencia e Innovaci\'on, Spain, 
projects MTM2015-67142-P (MINECO/FEDER, UE) and PGC2018-098279-B-I00 (MCIU/AEI/FEDER, UE). 
NMT thanks CWI, Amsterdam, for scientific support.

\end{document}